\documentclass[a4paper,12pt, reqno]{amsart}
\usepackage{a4wide}
\usepackage{amsthm}
\usepackage{amssymb}
 \newtheorem{theorem}{Theorem}
 
 \newtheorem{lemma}{Lemma}

 \newcommand{\mbb}{\mathbb}

 \newcommand{\Z}{\mbb{Z}}

 \theoremstyle{remark}

 \numberwithin{equation}{section}
\def\di{{\triangle}}   
\def\e{{\mathrm{e}}}

 \newcommand{\E}{\zeta_{Q} (s)} 
\begin{document}

%
%
%

\title{A note on the gaps between zeros of Epstein's zeta-functions on the critical line}
\author{Stephan Baier, Srinivas Kotyada and Usha Keshav Sangale}
\address{Jawaharlal Nehru University Delhi, School of Physical Sciences, Delhi 110067, India}
\address{Institute of Mathematical Sciences, HBNI, 
CIT Campus, Taramani, Chennai 600 113, India}
\address{SRTM University, Nanded, Maharashtra 431606, India}
\email[Stephan Baier]{email\textunderscore baier@yahoo.de}
\email[Kotyada Srinivas]{srini@imsc.res.in}
\email[Usha Sangale]{ushas073@gmail.com}
\begin{abstract} It is proved that Epstein's zeta-function $\zeta_{Q}(s)$, related to a positive definite integral binary quadratic form, has a zero $1/2 + i\gamma$ with $ T \leq \gamma \leq T + T^{{3/7} +\varepsilon} $ for sufficiently large positive numbers $T$. This is an improvement of the result by M. Jutila and K. Srinivas (Bull. London Math. Soc. 37 (2005) 45--53).
 \end{abstract}

 \maketitle
\hfill \textit{To Professor Matti Jutila with deep regards}

\section{ Introduction}

Let the quadratic form $ Q(x,y) = ax^2 + bxy + cy^2$ be positive definite, have integer coefficients and let $ r_Q(n)$ count the number of solutions of the equation $ Q(x,y) = n $ in integers $x$ and $y$. The Epstein zeta-function associated to $Q$ is denoted by $\zeta_Q(s)$ and is given by the series 
$$  \zeta_Q(s) = \sum_{(x,y)\in {\Z}^2 - (0,0)}\frac{1}{{Q(x,y)}^s} = \sum_{n=1}^{\infty}\frac{r_Q(n)}{n^s}
$$
in the half-plane $ \sigma > 1$, where (as usual) $ s=\sigma + it$. Throughout the paper, we shall write $\Delta$ to denote the number $ \Delta := \vert 4 ac - b^2 \vert$, the modulus of the discriminant of $Q$. $\zeta_Q(s)$ has many analytical properties in common with the Riemann zeta-function, $\zeta(s)$. For example, it admits analytic continuation into the entire complex plane except for a simple pole at $ s=1$ with residue $ 2\pi {\Delta}^{-1/2}$. It satisfies the following functional equation
\begin{equation}\label{f-eqn}
{\left( \frac {\sqrt{\Delta}}{2\pi}    \right)}^{s} \Gamma (s) \zeta_Q(s) = {\left( \frac {\sqrt{\Delta}}{2\pi}    \right)}^{1-s} \Gamma (1-s) \zeta_Q(1-s).
\end{equation}

\noindent
The analogue of Hardy's theorem for $\zeta(s)$ also holds true for $\E$, i.e., $\E$ admits infinitely many zeros on the critical line $\sigma = 1/2$. In fact, much more is true. In 1934, Potter and Tichmarsh \cite{P-T} showed that every interval of the type $ [T, T + T^{1/2 + \varepsilon}] $ contains a zero $ 1/2 + i\gamma$ of $\E$ for any fixed $\varepsilon$ and for all sufficiently large $T$. Sankaranarayanan in 1995 \cite{sank} showed that the same result holds true for intervals of the type $[ T, T+c T^{1/2}\log T]$. In 2005, Jutila and Srinivas \cite{J-S} proved that the same is true for intervals of the type $ [T,T + c T^{5/11 +\varepsilon}] $, thus surpassing the classical barrier of $1/2$ in exponent of $T$. In this note we improve this result further. More precisely, we prove the following 
\begin{theorem}
Let $Q$ be a positive definite binary integral quadratic form. Then for any fixed $\varepsilon > 0$ and $T\geq T(\varepsilon,Q)$, there is a zero $1/2+i\gamma$ of the corresponding Epstein zeta-function $\E$ with
\begin{equation}
\mid \gamma - T \mid \leq T^{3/7 + \varepsilon}.
\end{equation}
\end{theorem}

In this paper we  indicate  those steps in \cite{J-S} which enabled us to improve the result of Jutila and Srinivas mentioned earlier. Therefore, for technical details the readers are urged to refer to \cite{J-S} and \cite{J}. However, for the sake of completeness, we shall discuss the main ideas contained in \cite{J-S}. 
\medskip

The paper is organized as follows: In section \ref{basic-proof} we describe the basic idea of the proof, section \ref{lemmas} contains basic results used in the proof of the main theorem and in section \ref{d-sum} we estimate a double exponential sum non-trivially, which leads to the improvement.

\section{Basic idea of the proof}\label{basic-proof}

Hardy and Littlewood \cite{H-L} developed a beautiful method to prove the existence of a zero of the Riemann zeta-function $\zeta(s)$ on the critical line in a short interval. The significance of their method is that it is amenable to generalization.  We start by defining the functions $f(s), \gamma (s)$ and $W(t)$ as
$$
f(s)=\mathrm{e}^{\frac{1}{2} \pi i (\frac{1}{2}-s)}\left(
\frac{\sqrt{\di}}{2\pi}\right) ^s \Gamma(s) \E = \gamma (s) \zeta_Q(s)
$$
and
$$
W(t)=f\left(\frac{1}{2} +it\right).
$$
From the functional equation \eqref{f-eqn}, it follows that $W(t)$ is real for real values of $t$. Thus, the real zeros of $W(t)$ coincide with the zeros of $\zeta_Q(s)$ on the critical line. 

\medskip

First, let us assume that $W(t)$ has no zero in the interval $[ T-H, T+H ],$ with $T^{3\varepsilon} \leq H \leq T^{1/2}.$  Let $ H_0 = H
T^{-\epsilon}$ and consider the integral
$$
I= \int_{-H} ^{H} W(T+u) \e^{-(u/H_0 )^2} du .
$$
Then by our assumption 
\begin{equation}\label{equality}
|I|=\int ^H _{-H} |W(T+u)|\e^{-(u/H_0 )^2} du .
\end{equation}
If the equality in \eqref{equality} is violated, this will establish the existence of an odd order zero of $W(t)$ in the interval 
$[T-H, T+H]$. This contradiction is achieved by  estimating the integral in \eqref{equality} from below and above, provided $H = T^{3/7 + \varepsilon}$.

\medskip

Estimation from below is the easy step, thanks to a general result of K. Ramachandra \cite{kram} which states that the first 
power mean of a generalized Dirichlet series satisfying certain conditions can not be too small. We need only a particular case of this theorem
which is readily available as  Theorem 3 of \cite{B}, which we state as:

\begin{lemma}\label{kram}
Let $  B(s) = \sum_{n=1}^{\infty} b_n n^{-s} $ be any Dirichlet series satisfying the following conditions:
\begin{itemize}
\item[(i)] not all $b_n$'s are zero;
\item[(ii)] the function can be continued analytically in $ \sigma \geq a, \  |t| \geq t_0 $, and in this
region $ B(s) = O( (|t|+10)^A ).$
\end{itemize}
Then for every $\epsilon > 0, $ we have 
$$  \int_T^{T+H} | B(\sigma + it)| dt \gg H $$
for all $ H \geq (\log T)^{\epsilon}, \  T\geq T_0(\epsilon),$ and $\sigma > a $.
\end{lemma}
Thus the lower bound
\begin{equation}\label{the1e}
|I|\gg H_0 
\end{equation}
follows directly from Lemma \ref{kram}.

\medskip

Estimation from above is the hard part. We start with writing the integral $I$ as
\begin{equation}\label{the2e}
\begin{gathered}
I=\int_{-H} ^H \e^{\tfrac{1}{2} \pi
(T+u)}\left(\tfrac{\sqrt{\di}}{2\pi}\right)^{\frac{1}{2} +
i(T+u)}\Gamma\left(\tfrac{1}{2} +i(T+u)\right)
\\ \times \,\zeta_{Q}\left(\tfrac{1}{2} + i(T+u)\right)\e^{-(u/H_0)^2} du.
\end{gathered}
\end{equation}
The zeta-function $\zeta_{Q}\left(\tfrac{1}{2} + i(T+u)\right)$ appearing in the integrand is now replaced with an appropriate approximate formula. Such a formula was derived in \cite{J-S}, Lemma 1. We state this as 
\begin{lemma}\label{approx}
Let $t\geq2$ and $t^2\ll X \ll t ^A $, where $A$ is an arbitrarily
large positive constant. Then we have
\begin{equation}\label{the3e}
\begin{gathered}
\zeta_{Q}\left(\tfrac{1}{2} + it\right)=\sum_{n \leq X} r_Q (n)
n^{-1/2 -it}
\\ +\,(\log2)^{-1}\sum_{X<n\leq2X} r_Q(n) \log (2X/n)
n^{-1/2 - it}\\ +\,(\log 2)^{-1} 2\pi \di^{-1/2}
\left(\tfrac{1}{2} -it\right)^{-2}((2X)^{1/2 -it} -X ^{1/2 -it})+
O(tX^{-1/2}).
\end{gathered}
\end{equation}
\end{lemma}
Putting $X=T^3$, we observe that for $t\asymp T$ the last two error terms in \eqref{the3e} are $O(T^{-1/2})$. Substituting the  approximate formula for $\zeta_{Q}\left(\tfrac{1}{2} + i(T+u)\right)$ in  \eqref{the2e}  we, therefore,  have
%
\begin{eqnarray}\label{pf5e}
&& I = \sum_{{n\le T^3}}r_Q(n) n^{-1/2-iT}
             \int_{-H}^{H}\gamma(1/2+i(T+u))n^{-iu}e^{-(u/H_0)^2}du \nonumber\\
&&      +       (\log 2)^{-1}\sum_{T^3<n\le 2T^3}r_Q(n)\log (2T^3/n)n^{-1/2-it} \nonumber \\
&&     \times  \int_{-H}^{H}\gamma(1/2+i(T+u))n^{-iu}e^{-(u/H_0)^2}du + O(1).
\end{eqnarray}
As in \cite{J-S},   we extract, from the right hand side,  a weighted sum  of the form
\begin{equation}\label{int3e}
\sum_{n} \eta(n) r_{Q} (n) n^{-1/2 -it},
\end{equation}
where the weight function $\eta(n)$ is supported in the interval $[\,T\sqrt{\triangle}/2\pi - K, \, T\sqrt{\triangle}/2\pi + K \,]$, $t \asymp T$
and $t$ lies close to $T$.  The object is to show that this sum is \emph{small} in a certain sense. The remaining terms 
are evaluated by complex integration technique and we shall show that their contribution is negligible. 

\medskip

To begin with, the smooth weight function is $\eta(n)$ is defined as 
\[
\eta (x) =\left\{\begin{array}{ll} 
1 & \, \mathrm{for} \qquad | x-T\sqrt{\di}/2\pi |\leq K/2, \\
0 & \, \mathrm{for} \qquad |x-T\sqrt{\di}/2\pi|\geq K 
\end{array} \right.
\]
and $K$ is chosen to satisfy the relation 
\begin{equation}\label{the4e}
HK=T^{1+ 2\varepsilon}.
\end{equation}
Thus, using the trivial identity  $1 =  \eta(n) + (1 - \eta (n))  $ in \eqref{pf5e}, we obtain
\begin{eqnarray}\label{pf6e}
&& I = \sum_{\substack{n\le T^3 \\ |n-T\sqrt{\di}/2\pi|>K/2}}r_Q(n) ( 1 - \eta (n)) n^{-1/2-iT}
             \int_{-H}^{H}\gamma(1/2+i(T+u))n^{-iu}e^{-(u/H_0)^2}du \nonumber \\
&&      +      \int_{-H}^{H}\gamma(1/2+i(T+u))
       \left(\sum_{|n-T\sqrt{\di}/2\pi | \le K}r_Q(n) \eta(n) n^{-1/2-i(T+u)}\right) e^{-(u/H_0)^2}du \nonumber \\
&&      +       (\log 2)^{-1}\sum_{T^3<n\le 2T^3}r_Q(n)\log (2T^3/n)n^{-1/2-iT} \nonumber \\
&&     \times  \int_{-H}^{H}\gamma(1/2+i(T+u))n^{-iu}e^{-(u/H_0)^2}du + O(1) \nonumber \\
&&      =       S_1+S_2+S_3+O(1).
\end{eqnarray}
Now, we break $S_1$ into two sub-sums and estimate each of them individually.
\begin{eqnarray*}\label{pf7e}
&& S_1 = \sum_{T\sqrt{\di}/2\pi + K/2 \le n \le T^3}r_Q(n) (1 - \eta (n)) n^{-1/2-iT}
             \int_{-H}^{H}\gamma(1/2+i(T+u))n^{-iu}e^{-(u/H_0)^2}du\\
&&  +  \sum_{1 \le n \le T\sqrt{\di}/2\pi -  K/2}r_Q(n) (1 - \eta (n)) n^{-1/2-iT}
             \int_{-H}^{H}\gamma(1/2+i(T+u))n^{-iu}e^{-(u/H_0)^2}du\\
&&           = S_{11} + S_{12}.
\end{eqnarray*}

\noindent
\textbf{Estimation of $S_{11}$.}  The integral 
\begin{equation}\label{int-11}
 \int_{-H}^{H}\gamma(1/2+i(T+u))n^{-iu}e^{-(u/H_0)^2}du
\end{equation}
is first written as a contour integral on the rectangle with vertices $ \pm H, \pm H-iH_0$. 
Now using the well-know Stirling's formula for $\Gamma (s)$ which states that 
in any fixed vertical strip 
$-\infty < \alpha \leq \sigma \leq \beta < \infty$,
$$\Gamma(\sigma+it)=(2\pi)^{1/2} t^{\sigma+it-1/2}e^{-\frac{\pi}{2} t + 
\frac{\pi}{2} i (\sigma-1/2)-i t} ( 1 + O(1/t) ) \quad \hbox{as} \ t \rightarrow \infty ,$$
we obtain,
\begin{equation}\label{gamma}
\gamma(1/2+i(T+u)) n^{-iu}=\Delta^{1/4}\exp \left\lbrace 
i\left( T\log\left( {T\sqrt{\Delta}}/{2\pi}\right) -T+u\log \left( {T\sqrt{\Delta}}/{2\pi n}\right) 
     +O(1)\right) \right\rbrace . 
\end{equation}
On the vertical line, we have $u = -H + iy,  -H_0 \leq y \leq 0;$ and therefore,
\[
e^{-(u/H_0)^2} = e^{{- (H^2-y^2) + 2 i y H}/{(H_0)^2}} \ll e^{- (H/H_0)^2} \ll e^{-T^{2\varepsilon}}
\]
where as \eqref{gamma} is bounded. On the other hand on the lower horizontal side of the rectangle, setting
$ u = x - i H_0, -H \leq x \leq H$, we observe that 
\[
e^{-(u/H_0)^2} \leq e^{{- (x^2-(H_0)^2)/(H_0)^2}} \leq e^{(H_0)^2/(H_0)^2} = 1
\]
and
\begin{equation}\label{gamma-2}
\gamma(1/2+i(T+u)) n^{-iu} \ll e^{H_0 \log \left( {T\sqrt{\Delta}}/{2\pi n}\right)} = e^{{- H_0 \log \left( {{2\pi n}/ T\sqrt{\Delta}} \right)}}.
\end{equation}
For $ n \geq T\sqrt{\Delta}/{2\pi} + K/2$, using the elementary inequality $\log a/b \geq |a-b|/(a+b)$, we see that
$$H_0 \log \left( {2\pi n}/ T\sqrt{\Delta} \right) \geq H_0K/ 2(T\sqrt{\Delta}/{2\pi} + K/2).$$
Now, 
\[
H_0K/2(T\sqrt{\Delta}/{2\pi} + K/2) \geq T^{\varepsilon}  \quad  \textrm{provided} \quad H_0K > T^{1 + \varepsilon}.
\]
The last inequality is guaranteed by the choice of $K$ in \eqref{the4e}. Therefore, from \eqref{gamma-2}, we get
\[
\gamma(1/2+i(T+u)) n^{-iu} \ll e^{-T^{\varepsilon}}
\]
Collecting all the estimates above, we have for $ n \geq T\sqrt{\Delta}/{2\pi} + K/2$,
\[
 \int_{-H}^{H}\gamma(1/2+i(T+u))n^{-iu}e^{-(u/H_0)^2}du \ll H e^{-T^{\varepsilon}}.
\]
Thus,
\begin{equation}\label{est-11}
\begin{gathered}
S_{11} \ll  H e^{-T^{\varepsilon}} \vert \sum_{T\sqrt{\di}/2\pi + K/2 \le n \le T^3}r_Q(n) (1 -\eta(n) )n^{-1/2-iT} \vert \\
\ll H e^{-T^{\varepsilon}}  \sum_{T\sqrt{\di}/2\pi + K/2 \le n \le T^3}  r_Q(n)  n^{-1/2} \\
\ll H e^{-T^{\varepsilon}}T^{3+\epsilon} (T\sqrt{\di}/2\pi + K/2)^{-1/2} \ll 1.
\end{gathered}
\end{equation}

\noindent
\textbf{Estimation of $S_{12}$.} The estimation of $S_{12}$ is similar to that of $S_{11}$. In this case, the integral \eqref{int-11} is written as a contour integral over a rectangle with vertices $ \pm H, \pm H +i H_0$. Then the integral is bounded by $He^{-T^{\varepsilon}}$ on the vertical sides, where as on the horizontal sides, putting $ u = x + iH_0, -H \leq x \leq H;$ we  get
$$
| \gamma(1/2+i(T+u)) n^{-iu} | \leq e^{-H_0 \log \left( {T\sqrt{\Delta}}/{2\pi n}\right)} \ll e^{-T^{\varepsilon}}.
$$
Therefore,
\begin{equation}
\begin{gathered}\label{est-12}
| S_{12} | \ll H e^{-T^{\varepsilon}} \vert \sum_{1 \le n \le T\sqrt{\di}/2\pi -  K/2}r_Q(n) (1 -\eta(n) ) n^{-1/2-iT} \vert \\
\ll H e^{-T^{\varepsilon}} \sum_{1 \le n \le T\sqrt{\di}/2\pi -  K/2}  r_Q(n) n^{-1/2} \\
\ll H e^{-T^{\varepsilon}} T^{1/2 +\varepsilon} \\
 \ll 1.
\end{gathered}
\end{equation}
 \textbf{Estimation of $S_3$.} The estimation of $S_3$ follows the same pattern as that of $S_{11}$. To show that the integral \eqref{int-11} is small, we want 
 $$ 
 H_0 \log \left( {2\pi n}/ T\sqrt{\Delta} \right) > T^{\epsilon},
 $$
 that is
 $$
 H_0 \log \left( \tfrac{{2\pi n} -T\sqrt{\Delta}}{ {{2\pi n} +T\sqrt{\Delta}}} \right) > T^{\epsilon},
 $$
 or 
 $$
 H_0 \log \left( \tfrac{{2\pi T^3} -T\sqrt{\Delta}}{ {{2\pi T^3} +T\sqrt{\Delta}}} \right) > T^{\epsilon},
 $$
 which is the same as
 $$
 H_0 \log \left( 1 - \tfrac{2T\sqrt{\Delta}}{ {2\pi T^3} +T\sqrt{\Delta}} \right) > T^{\epsilon},
 $$
 which is true if $H_0 > T^{2\varepsilon}$. This is indeed the case, since our choice of $H$ is $ T^{3\varepsilon} \leq H \leq T^{1/2}$. Thus we conclude that
 \begin{eqnarray}\label{est-s3}
| S_3 |  & \ll & (\log 2)^{-1}\sum_{T^3<n\le 2T^3} T^{\varepsilon} (\log 2) \, T^{-3/2}  | \int_{-H}^{H}\gamma(1/2+i(T+u))n^{-iu}e^{-(u/H_0)^2} | du \nonumber \\
& \ll & H e^{-T^{\epsilon}} T^{\varepsilon}T^3 T^{-3/2} \ll 1.
\end{eqnarray}

Now it remains to estimate $S_2$.

\medskip

\noindent
\textbf{Estimation of $S_2$.}
We have
\begin{eqnarray*}
S_2 & = &\int_{-H}^{H}\gamma(1/2+i(T+u))
       \left(\sum_{|n-T\sqrt{\di}/2\pi | \le K}r_Q(n)\eta(n) n^{-1/2-i(T+u)}\right) e^{-(u/H_0)^2}du \\
       &= &\int_{T-H}^{T+H}\gamma(1/2+it)
       \left(\sum_{|n-T\sqrt{\di}/2\pi | \le K}r_Q(n)\eta(n) n^{-1/2-it}\right) e^{-((t-T)/H_0)^2}dt \\
       &= &\int_{T-H}^{T+H}\gamma(1/2+it)
       \left(\sum_{n=1}^{\infty}r_Q(n)\eta(n)n^{-1/2-it}\right) e^{-((t-T)/H_0)^2}dt. \\
\end{eqnarray*}
Therefore,
\begin{equation}\label{est-s2}
|S_2| \leq H \sup_{|T-t| \leq H} \mid \sum_{n=1}^{\infty}r_Q(n)\eta(n)n^{-1/2-it} \mid
\end{equation}
The objective now is to show that 
\begin{equation}\label{smooth-sum}
\sum_{n=1}^{\infty}r_Q(n)\eta(n)n^{-1/2-it} \ll (\log T)^{-2}
\end{equation}
for a suitable choice of the parameter $K$. Then combining this with \eqref{pf6e}, \eqref{est-11}, \eqref{est-12}, \eqref{est-s3} and \eqref{est-s2}, we have
\[
| I | \ll H_0 (\log T)^{-2}.
\]
This is a contradiction to \eqref{the1e}. 

\medskip
In \cite{J-S}, the crucial sum in \eqref{smooth-sum} 
 was first transformed into another sum (equation (3.6) of \cite{J-S}) using a transformation formula. By partial summation, this sum got reduced to the estimation of the following expression (for notations see subsection 4.1):
\begin{equation}\label{crucial-2}
K^{1/4}N^{-1/4}T^{-1/2} \sum_{N\leq Q^{\ast}(x,y) \leq N'} e\left(Q^{\ast}(x,y) \left( \frac{\overline{h\Delta_0}}{k}-\frac{1}{2hk\Delta_0} \right)+\frac{t}{\pi}\cdot \phi\left(\frac{\pi Q^{\ast}(x,y)}{2hk\Delta_0 t}\right)\right)
\end{equation}
This is  equation (4.4) of \cite{J-S}. The task is show that this term is $O(T^{-\varepsilon})$. To establish this, the double exponential sum, appearing above, was estimated non-trivially in one variable using van der Corput's method and trivial estimate was taken in the other variable. The authors obtained the bound $O(K^{11/12}T^{-1/2})$, which is $O(T^{-\varepsilon})$, provided $K= T^{6/11 - \varepsilon}$.

In the present paper we estimate the above double exponential sum  non-trivially in both variables. To show that \eqref{crucial-2} is $O(T^{-\varepsilon})$, it is now enough to take $K=T^{4/7 - \varepsilon}$ and thus the gap $H = T^{3/7 + \varepsilon}$ follows.

\section{Preliminary lemmas}\label{lemmas}

We will use the following well-known lemmas in the proof of our theorem.

\begin{lemma} (Generalized Weyl differencing)\label{W}
 Let $a < b$ be integers, $\lambda$ be a natural number and  $\xi(n)$ be a complex valued function such that $\xi(n)=0$ if $n\not\in (a,b]$. If $H$ is a positive integer then
\begin{equation}\label{Weyl-difference}
{\mid \sum_n \xi(n) \mid}^2 \leq \frac{(b-a) + H}{H} \sum_{|h|<H} \left( 1 - \frac{\lambda |h|}{H} \right) \sum_n \xi(n)\overline{\xi (n-\lambda h)}
\end{equation}
\end{lemma}

\noindent
\textbf{Proof.} For the case  $\lambda =1$, this is Lemma 2.5 of \cite{graham}. The general case can be proved similarly.

\begin{lemma} (B-Process , Lemma 3.6 of \cite{graham} )\label{B}
Suppose that $f$ has four continuous derivatives on $[a, b]$, and that $f''<0$ on this interval. Suppose further that $[a, b] \subseteq [N, 2N]$ and that $ \alpha = f'(b)$ and $\beta = f'(a)$. Assume that there is some $F>0$ such that
\begin{equation*}
f^{(2)}(x) \asymp FN^{-2}, \ f^{(3)}(x) \ll F N^{-3}, \textrm{and} \  f^{(4)}(x) \ll F N^{-4}
\end{equation*}
for $x$ in $[a, b]$. Let $x_{\nu}$ be defined by the relation $ f'(x_{\nu}) = \nu$, and let $\phi(\nu) = - f(x_{\nu}) + \nu x_{\nu}$. Then
\begin{equation}\label{b-process}
\sum_{a\leq n \leq b} e(f(n)) = \sum_{\alpha \leq \nu \leq \beta} \frac{ e(-\phi(\nu) - 1/8 )}{ {|f''(x_{\nu}|}^{1/2} } + O(\log(FN^{-1} +2 ) + F^{-1/2}N ).
\end{equation}
\end{lemma}

\begin{lemma} (van der Corput's bound, Theorem 2.2 of \cite{graham})\label{C}
Suppose that $f$ is a real valued function with two continuous derivatives on $[a, b]$ where $ a<b$ are integers. Suppose also that there is some $\lambda > 0$ and some $ \alpha \geq 1$ such that
$$
\lambda \leq | f''(x)| \leq \alpha \lambda
$$
on $[a, b]$. Then
\begin{equation}\label{van der Corput}
\sum_{a\leq n \leq b} e(f(n)) \ll \alpha (b-a) {\lambda}^{1/2} + {\lambda}^{-1/2}.
\end{equation} 
\end{lemma}

\section{Estimation of the double exponential sum}\label{d-sum}

\subsection{Preparation and description of the method}
First, we explain the meanings of the functions and variables occurring in \eqref{crucial-2}. The function $Q^{\ast} (x,y)$ denotes a certain positive definite quadratic form (for the details of its definition, see equation (3.4) of \cite{J-S}) 
\begin{equation}
Q^{\ast}(x,y)=a^{\ast}x^2+b^{\ast}xy+c^{\ast}y^2, \quad a^{\ast},b^{\ast},c^{\ast}\in \mathbb{Z},
\end{equation} 
which is related to $Q(x,y)$ and a positive integer $k$ in a specific way. Throughout the sequel,  we denote by   $d$ the discriminant of this form, i.e. 
$$ d :=\left(b^{\ast}\right)^2-4a^{\ast}c^{\ast}<0. $$ 
The form $Q^{\ast}(x,y)$  is defined in such a way that $a^{\ast} >0, c^{\ast}>0$. Further, as remarked in \cite{J-S},  $ | d | \leq \Delta. $
 
 \medskip
 
 \noindent
The variables in \eqref{crucial-2} satisfy the following conditions. We suppose that $1\leq N \leq N' \leq 2N$, $N\ll K$, $\Delta_0$, $h$ and $k$ are positive integers 
satisfying $\Delta_0|\Delta$, $ (h \Delta_0,k) =1$ and 
$$
\frac{K}{T} \ll \mid \frac{1}{\sqrt{\Delta}} - \frac{h}{k} \mid \leq \frac{\pi K}{T\Delta},
$$ 
and the sizes of $k,h$ and the real number $t$ are 
$$k, h \asymp \sqrt{T/K}, \quad t \asymp T.$$ 
As usual, $\overline{h\Delta_0}$ denotes a multiplicative inverse of $h\Delta_0$ modulo $k$. 

\medskip

\noindent
Finally, the function $\phi(x)$ is defined as
$$
\phi(x) = \mbox{arsinh} (x^{1/2}) + (x + x^2 )^{1/2}.
$$

\medskip

\noindent
Our goal is now to bound non-trivially the exponential sum in \eqref{crucial-2}, i.e. the exponential sum
\begin{equation}\label{main-sum}
\sum\limits_{x} \sum\limits_{y\in \mathcal{I}(x)} e\left(Q^{\ast}(x,y)\cdot r+\frac{t}{\pi}\cdot \phi\left(\frac{\pi Q^{\ast}(x,y)}{2hk\Delta_0 t}\right)\right),
\end{equation}
where 
$$
r:=\frac{\overline{h\Delta_0}}{k}-\frac{1}{2hk\Delta_0}
$$
and 
$$
\mathcal{I}(x):=\{y\in \mathbb{R} \ :\ N\le Q^{\ast}(x,y)\le N'\}.
$$ 
In \cite{J-S}, the summation over $y$ was evaluated using the following classical  estimate for exponential sums (Lemma 4.1 of \cite{J-S}). 

\begin{lemma}
Suppose that $f$ is a real valued function with three continuous derivatives on $[a, b]$ where $ a<b$ are integers. 
Suppose also that there is some $\lambda > 0$ and some $ \alpha \geq 1$ such that
$$
\lambda \leq | f'''(x)| \leq \alpha \lambda
$$
on $[a, b]$. Then
$$
\sum_{a\leq n \leq b} e(f(n)) \ll \alpha^{1/2} (b-a) {\lambda}^{1/6} + (b-a)^{1/2}{\lambda}^{-1/6}.
$$
\end{lemma}

The above lemma can be proved using Weyl differencing, Lemma \ref{W}, followed by applying the can der Corput bound, Lemma \ref{C}. 
\medskip

\noindent
In \cite{J-S},
the sum over $x$ was treated trivially. In the present paper, we also want to exploit cancellations in the $x$-sum. To this end,
we explicitly carry out Weyl differencing for the sum over $y$, then employ the B process, Lemma \ref{B}, re-arrange the summation and finally
apply van der Corput's bound to the sum over $x$. 

\medskip

\noindent
It is easy to see that $\mathcal{I}(x)$ is empty unless $x\in \mathcal{J}$, where
$$
\mathcal{J}:=\left[-\frac{2\sqrt{c^{\ast}N'}}{\sqrt{|d|}},\frac{2\sqrt{c^{\ast}N'}}{\sqrt{|d|}}\right]
$$
and that 
$$
\mathcal{I}(x)=I(x) \cup I'(x),
$$
where 
$$
I(x):=\left[\sqrt{\max\left\{0, \frac{N}{c^{\ast}}-\frac{|d|}{\left(2c^{\ast}\right)^2}\cdot x^2\right\}}, 
\sqrt{\frac{N'}{c^{\ast}}-\frac{|d|}{\left(2c^{\ast}\right)^2}\cdot x^2}\right] 
$$
and 
$$
I'(x):=\left[-\sqrt{\frac{N'}{c^{\ast}}-\frac{|d|}{\left(2c^{\ast}\right)^2}\cdot x^2},-\sqrt{\max\left\{0,\frac{N}{c^{\ast}}-\frac{|d|}{\left(2c^{\ast}\right)^2}\cdot x^2\right\}}
\right].
$$

Set 
$$
J=\left\{x\in \mathcal{J} \ :\ x\ge 0\right\}=\left[0,\frac{2\sqrt{c^{\ast}N'}}{\sqrt{|d|}}\right].
$$ 
In the following, we estimate the partial sum
\begin{equation*}
\sum\limits_{x\in J} \sum\limits_{y\in I(x)} e\left(Q^{\ast}(x,y)\cdot r+\frac{t}{\pi}\cdot \phi\left(\frac{\pi Q^{\ast}(x,y)}{2hk\Delta_0 t}\right)\right),
\end{equation*}
where $x$ is positive and $y$ runs over the interval $I(x)$. The remaining three partial sums with a) $x\ge 0$ and $y\in I'(x)$, b) $x<0$ and $y\in I(x)$, c) $x<0$ and $y\in I'(x)$ can be estimated in a similar way. 

\subsection{Application of Weyl differencing}
We start with applying the Cauchy-Schwarz inequality, getting 
\begin{equation} \label{cs1}
\begin{split}
& \left|\sum\limits_{x\in J} \sum\limits_{y\in I(x)} e\left(Q^{\ast}(x,y)\cdot r+\frac{t}{\pi}\cdot \phi\left(\frac{\pi Q^{\ast}(x,y)}{2hk\Delta_0 t}\right)\right)\right|^2\\
\ll & \sqrt{N}\cdot 
\sum\limits_{x\in J} \left|\sum\limits_{y\in I(x)} e\left(Q^{\ast}(x,y)\cdot r+\frac{t}{\pi}\cdot \phi\left(\frac{\pi Q^{\ast}(x,y)}{2hk\Delta_0 t}\right)\right)\right|^2,
\end{split}
\end{equation}
where we use $|J|\ll \sqrt{N}$. 
Applying Lemma \ref{W} with $\lambda = 2$, and using $|J|\ll \sqrt{N}$ and $| I(x)|\ll \sqrt{N}$, we have 
\begin{equation} \label{weyl}
\begin{split}
& \sum\limits_{x\in J}\left|\sum\limits_{y\in I(x)} e\left(Q^{\ast}(x,y)\cdot r+\frac{t}{\pi}\cdot \phi\left(\frac{\pi Q^{\ast}(x,y)}{2hk\Delta_0 t}\right)\right)\right|^2\\
\ll  & \sum\limits_{x\in J}\frac{\sqrt{N}}{M} \cdot \sum\limits_{0\le |m|\le M} \left(1-\frac{|2m|}{M}\right) \cdot \sum\limits_{y\in I_m(x)} e\left(f_x(y+2m)-f_x(y)\right)\\
\ll & \frac{\sqrt{N}}{M} \cdot \sum\limits_{1\le m\le M} \left| \sum\limits_{x\in J} \sum\limits_{y\in I(x)} e\left(f_x(y+m)-f_x(y-m)\right)\right| + NM+\frac{N^{3/2}}{M},
\end{split}
\end{equation} 
where $M\le N$ is any natural number, 
$$
I_m(x):=\{y\in I(x) \ :\ y+2m\in I(x)\}
$$
and 
$$
f_x(y):=Q^{\ast}(x,y)\cdot r+\frac{t}{\pi}\cdot \phi\left(\frac{\pi Q^{\ast}(x,y)}{2hk\Delta_0 t}\right).
$$

We have  
$$
Q^{\ast}(x,y+m)-Q^{\ast}(x,y-m)=2m\left(b^{\ast}x+2c^{\ast}y\right)
$$
and, using Taylor series expansion, 
\begin{equation*} 
\phi\left(\frac{\pi Q^{\ast}(x,y+m)}{2hk\Delta_0 t}\right)-\phi\left(\frac{\pi Q^{\ast}(x,y-m)}{2hk\Delta_0 t}\right)=2mg_x'(y)+O\left(m^3\left|g_x'''(y)\right|\right),
\end{equation*}
where 
$$
g_x(y):=\phi\left(\frac{\pi Q^{\ast}(x,y)}{2hk\Delta_0 t}\right).
$$
We note that the quadratic Taylor term disappears due to our treatment. This is the reason why we chose $\lambda=2$ in our application of Lemma \ref{W}. The
absence of a quadratic term is advantageous for us because the cubic term can be handled easily, whereas the presence of a quadratic term
would lead to difficulties. 

Using 
$$
\phi'(u)=\frac{1}{\sqrt{u}}+O(\sqrt{u}), \quad \phi''(u)=O\left(\frac{1}{u^{3/2}}\right) \quad \mbox{and} \quad \phi'''(u)=O\left(\frac{1}{u^{5/2}}\right)
$$
for $|u|\le 1$, we calculate that 
$$
g_x'(y)=\frac{\pi\left(b^{\ast}x+2c^{\ast}y\right)}{2hk\Delta_0 t}\cdot \left(\frac{1}{\sqrt{\frac{\pi Q^{\ast}(x,y)}{2hk\Delta_0 t}}} +O\left(\sqrt{\frac{\pi Q^{\ast}(x,y)}{2hk\Delta_0 t}}\right)\right)=\frac{\sqrt{\pi}\left(b^{\ast}x+2c^{\ast}y\right)}{\sqrt{2hk\Delta_0 tQ^{\ast}(x,y)}}+O\left(\frac{NK^{3/2}}{T^3}\right)
$$
and 
$$
g_x'''(y)=O\left(\frac{K^{1/2}}{NT}\right).
$$
It follows that 
\begin{equation} \label{appro}
\sum\limits_{y\in I(x)} e\left(f_x(y+m)-f_x(y-m)\right)\\
= \sum\limits_{y\in I(x)} e\left(2mF_x(y)\right)+
O\left(\frac{mN^{3/2}K^{3/2}}{T^2}+\frac{m^3K^{1/2}}{N^{1/2}}\right),
\end{equation}
where 
$$
F_x(y):=\left(b^{\ast}x+2c^{\ast}y\right)\cdot \left(r+
\frac{\sqrt{t}}{\sqrt{2\pi hk\Delta_0 Q^{\ast}(x,y)}}\right).
$$

\subsection{Application of the B process}
Now we want to employ the B process, Lemma \ref{B}, to transform the exponential sum on the right-hand side of \eqref{appro}. 
To this end, we calculate that
\begin{equation} \label{Fprime}
\begin{split}
F_x'(y)= & 2c^{\ast}r+
\frac{\sqrt{t}|d|x^2}{2\sqrt{2\pi hk\Delta_0} Q^{\ast}(x,y)^{3/2}},\\
F_x''(y)= & -\frac{3\sqrt{t}|d|x^2\left(b^{\ast}x+2c^{\ast}y\right)}{4\sqrt{2\pi hk\Delta_0} Q^{\ast}(x,y)^{5/2}}\asymp \frac{K^{1/2}}{N},\\
F_x'''(y)= & O\left(\frac{K^{1/2}}{N^{3/2}}\right),\\
F_x''''(y)= & O\left(\frac{K^{1/2}}{N^{2}}\right).
\end{split}
\end{equation}
We also need to find the precise range in which $F_x'(y)$ lies, which we do in the following.
By our assumptions on $x$ and $y$, we have $N\le Q^{\ast}(x,y)\le N'$. For fixed $x$, we also have
$$
Q^{\ast}(x,y)\ge \frac{|d|}{4c^{\ast}}\cdot x^2. 
$$
Hence $a(x)\le F'_x(y)\le b(x)$, where  
$$
a(x):=2c^{\ast}r+
\frac{\sqrt{t}|d|x^2}{2\sqrt{2\pi hk\Delta_0} N'^{3/2}}
$$
and 
\begin{equation*}
\begin{split}
\quad b(x):= & 2c^{\ast}r+
\frac{\sqrt{t}|d|x^2}{2\sqrt{2\pi hk\Delta_0} \max\{N,|d|x^2/(4c^{\ast})\}^{3/2}}\\ = & 2c^{\ast}r+
\min\left\{\frac{\sqrt{t}|d|x^2}{2\sqrt{2\pi hk\Delta_0} N^{3/2}},\frac{\left(4c^{\ast}\right)^{3/2}\sqrt{t}}{2\sqrt{2\pi |d|hk\Delta_0} x}\right\}.
\end{split}
\end{equation*}
Now Lemma \ref{B} yields
\begin{equation} \label{expsums}
\begin{split}
& \sum\limits_{y\in I(x)} e\left(2mF_x(y)\right)\\
= & \sum\limits_{2ma(x)\le n\le 2mb(x)} \frac{e\left(2mF_x\left(y_{x,m,n}\right)-ny_{x,m,n}-1/8\right)}{\sqrt{2m|F_x''(y_{x,m,n})|}}+O\left(\log T+ \frac{N^{1/2}}{m^{1/2}K^{1/4}}\right),
\end{split}
\end{equation}
where $y_{x,m,n}\in I(x)$ is the solution of $2mF_x'(y_{x,m,n})=n$. 

We compute that
$$
F_x\left(y_{x,m,n}\right)=\sqrt{-|d|x^2+4c^{\ast}\left(\frac{m\sqrt{t}|d|x^2}{\sqrt{2\pi hk\Delta_0}\left(n-4mc^{\ast}r\right)}\right)^{2/3}}\cdot \left(r+\left(\frac{\left(n-4mc^{\ast}r\right)t}{2\pi hk\Delta_0|d|mx^2}\right)^{1/3}\right)
$$
and 
$$
ny_{x,m,n}=\frac{n}{2c^{\ast}}\cdot \left(-b^{\ast}x+\sqrt{-|d|x^2+4c^{\ast}\left(\frac{m\sqrt{t}|d|x^2}{\sqrt{2\pi hk\Delta_0}\left(n-4mc^{\ast}r\right)}\right)^{2/3}}\right).
$$
Putting together gives
\begin{equation} \label{putting}
\begin{split}
& G_{m,n}(x):=2mF_x\left(y_{x,m,n}\right)-ny_{x,m,n}\\
= & \frac{b^{\ast}n}{2c^{\ast}}\cdot x-\frac{1}{2c^{\ast}}\cdot \sqrt{-|d|x^2+4c^{\ast}\left(\frac{m\sqrt{t}|d|x^2}{\sqrt{2\pi hk\Delta_0}\left(n-4mc^{\ast}r\right)}\right)^{2/3}}\times\\ & 
\left(\left(n-4c^{\ast}mr\right)-4c^{\ast}m\left(\frac{\left(n-4mc^{\ast}r\right)t}{2\pi hk\Delta_0|d|mx^2}\right)^{1/3}\right)\\
= & \frac{b^{\ast}n}{2c^{\ast}}\cdot x+\frac{(n-4c^{\ast}mr)}{2c^{\ast}|d|x^2}\cdot \left(-|d|x^2+4c^{\ast}\left(\frac{m\sqrt{t}|d|x^2}{\sqrt{2\pi hk\Delta_0}\left(n-4mc^{\ast}r\right)}\right)^{2/3}\right)^{3/2}\\
= & \frac{b^{\ast}n}{2c^{\ast}}\cdot x+\frac{1}{2c^{\ast}}\cdot \left(-(n-4c^{\ast}mr)^{2/3}|d|^{1/3}x^{2/3}+4c^{\ast}\cdot \frac{m^{2/3}t^{1/3}}{(2\pi hk\Delta_0)^{1/3}}\right)^{3/2}.
\end{split}
\end{equation}

Combining  \eqref{cs1}, \eqref{weyl}, \eqref{appro}, \eqref{Fprime}, \eqref{expsums} and \eqref{putting}, we get
\begin{equation} \label{cs}
\begin{split}
& \left|\sum\limits_{x\in J} \sum\limits_{y\in I(x)} e\left(Q^{\ast}(x,y)\cdot r+\frac{t}{\pi}\cdot \phi\left(\frac{\pi Q^{\ast}(x,y)}{2hk\Delta_0 t}\right)\right)\right|^2\\
\ll &  \frac{N}{M} \cdot \sum\limits_{1\le m\le M} \left| \sum\limits_{x\in J} \sum\limits_{ma(x)\le n\le mb(x)} \frac{e\left(G_{m,n}(x)\right)}{\sqrt{|F_x''(y_{x,m,n})|}}\right|+\\ & 
O\left(MN^{3/2}+\frac{N^2}{M}+\frac{MN^{5/2}K^{3/2}}{T^2}+M^3N^{1/2}K^{1/2}+\frac{N^2}{M^{1/2}K^{1/4}}\right).
\end{split}
\end{equation}

\subsection{Application of van der Corput's bound}
Let's first work out what the trivial estimate for the double exponential sum above gives. We will see that we recover precisely 
the result
in \cite{J-S} in this way. Clearly, 
\begin{equation*}
|J|\ll N^{1/2} \quad \mbox{and} \quad b(x)-a(x)=O\left(K^{1/2}N^{-1/2}\right).
\end{equation*}
Together with \eqref{Fprime}, this implies
\begin{equation} \label{trivial}
 \frac{N}{M} \cdot \sum\limits_{1\le m\le M} \left| \sum\limits_{x\in J} \sum\limits_{2ma(x)\le n\le 2mb(x)} \frac{e\left(G_{m,n}(x)\right)}{\sqrt{2m|F_x''(y_{x,m,n})|}}\right| =O\left(M^{1/2}N^{3/2}K^{1/4}\right).
\end{equation}
Choosing $M:=\left[N^{1/3}K^{-1/6}\right]$ to balance the $O$-term in \eqref{trivial} and the second $O$-term $N^2/M$ in \eqref{cs}, 
using $N\ll K$, and taking the square root, we deduce that
\begin{equation}
\begin{split}
& K^{1/4}N^{-1/4}T^{-1/2}\left|\sum\limits_{x\in J} \sum\limits_{y\in I(x)} e\left(Q^{\ast}(x,y)\cdot r+\frac{t}{\pi}\cdot \phi\left(\frac{\pi Q^{\ast}(x,y)}{2hk\Delta_0 t}\right)\right)\right|\\
= & O\left(\frac{K^{11/12}}{T^{1/2}}+\frac{K^{25/12}}{T^{3/2}}\right).
\end{split}
\end{equation}
So if 
$$
K:=T^{6/11-\varepsilon},
$$
the above is $O\left(T^{-\varepsilon}\right)$, as in \cite{J-S}. 

Now we estimate the said double exponential sum non-trivially, thus obtaining an improvement over the result in \cite{J-S}. 
First, we re-arrange summations, getting 
\begin{equation} \label{rearrange}
\sum\limits_{x\in J} \sum\limits_{2ma(x)\le n\le 2mb(x)} \frac{e\left(G_{m,n}(x)\right)}{\sqrt{2m|F_x''(y_{x,m,n})|}}=\sum\limits_{n\in J_m} \sum\limits_{A_m(n)\le x\le B_m(n)}  \frac{e\left(G_{m,n}(x)\right)}{\sqrt{2m|F_x''(y_{x,m,n})|}},
\end{equation}
where
$$
J_m= \left[4c^{\ast}mr
,4c^{\ast}mr+\frac{2c^{\ast}m\sqrt{t}N'}{\sqrt{2\pi hk\Delta_0}N^{3/2}}\right],
$$
$$
A_m(n):=\left(n-4mc^{\ast}r\right)^{1/2}\cdot  \frac{(2\pi hk\Delta_0)^{1/4} N^{3/4}}{(m|d|)^{1/2}t^{1/4}}
$$
and
$$
B_m(n):=\min\left\{\frac{2(c^{\ast}N')^{1/2}}{|d|^{1/2}},\left(n-4mc^{\ast}r\right)^{1/2}\cdot  \frac{(2\pi hk\Delta_0)^{1/4} N'^{3/4}}{(m|d|)^{1/2}t^{1/4}},\frac{\left(4c^{\ast}\right)^{3/2}m t^{1/2}}{(2\pi |d|hk\Delta_0)^{1/2}\left(n-4mc^{\ast}r\right)}\right\}.
$$

Our idea is to use the van der Corput bound, Lemma \ref{C}, to estimate the inner sum over $x$ on the right-hand side of \eqref{rearrange}.
To this end, we note that
\begin{equation} \label{note}
|J_m|\ll \frac{mK^{1/2}}{N^{1/2}} \quad \mbox{and} \quad B_m(n)-A_m(n)\ll N^{1/2}
\end{equation}
and compute that
\begin{equation}
\begin{split}
G_{m,n}''(x)= & -\frac{(n-4c^{\ast}mr)^{2/3}|d|^{1/3}}{6c^{\ast}x^{4/3}}\cdot \left((n-4c^{\ast}mr)^{2/3}|d|^{1/3}x^{2/3}+4c^{\ast}\cdot \frac{m^{2/3}t^{1/3}}{(2\pi hk\Delta_0)^{1/3}}\right)\times\\ &  \left(-(n-4c^{\ast}mr)^{2/3}|d|^{1/3}x^{2/3}+4c^{\ast}\cdot \frac{m^{2/3}t^{1/3}}{(2\pi hk\Delta_0)^{1/3}}\right)^{-1/2}.
\end{split}
\end{equation}
Further, we observe that
\begin{equation*}
\frac{(n-4c^{\ast}mr)^{2/3}|d|^{1/3}}{6c^{\ast}x^{4/3}} \asymp \frac{m^{2/3}t^{1/3}}{(hk\Delta_0)^{1/3} N}\asymp \frac{m^{2/3}K^{1/3}}{N}
\end{equation*}
and 
\begin{equation*}
(n-4c^{\ast}mr)^{2/3}|d|^{1/3}x^{2/3}+4c^{\ast}\cdot \frac{m^{2/3}t^{1/3}}{(2\pi hk\Delta_0)^{1/3}} \asymp \frac{m^{2/3}t^{1/3}}{(hk\Delta_0)^{1/3}}\asymp m^{2/3}K^{1/3}
\end{equation*}
and hence
\begin{equation} \label{preest}
G_{m,n}''(x)\asymp \frac{m^{4/3}K^{2/3}}{N} \cdot \left(-(n-4c^{\ast}mr)^{2/3}|d|^{1/3}x^{2/3}+4c^{\ast}\cdot \frac{m^{2/3}t^{1/3}}{(2\pi hk\Delta_0)^{1/3}}\right)^{-1/2}
\end{equation} 
for $n$ and $x$ in the relevant summation intervals. 

We first assume that $n\in J'_m\subseteq J_m$, where 
\begin{equation*}
J'_m:=\left.\left[4c^{\ast}mr,4c^{\ast}mr+\frac{c^{\ast}m\sqrt{t}}{\sqrt{2\pi hk\Delta_0N'}}\right.\right),
\end{equation*}
in which case we compute that 
\begin{equation*}
-(n-4c^{\ast}mr)^{2/3}|d|^{1/3}x^{2/3}+4c^{\ast}\cdot \frac{m^{2/3}t^{1/3}}{(2\pi hk\Delta_0)^{1/3}} \asymp \frac{m^{2/3}t^{1/3}}{(hk\Delta_0)^{1/3}} \asymp m^{2/3}K^{1/3}
\end{equation*}
and hence, using \eqref{preest},
\begin{equation} \label{G''asymp}
G_{m,n}''(x) \asymp \frac{mK^{1/2}}{N} \quad \mbox{ if } A_m(n)\le x\le B_m(n).
\end{equation}
Now we apply partial summation to remove the factor $1/\sqrt{F_x''(y_{x,m,n})}$ and Lemma \ref{C} to get
\begin{equation} \label{aftervandercorput}
\begin{split}
& \sum\limits_{n\in J_m'} \sum\limits_{A_m(n)\le x\le B_m(n)}  \frac{e\left(G_{m,n}(x)\right)}{\sqrt{2m|F_x''(y_{x,m,n})|}} \\ \ll & \frac{N^{1/2}}{m^{1/2}K^{1/4}}\cdot \left(mK^{1/2}\cdot \left(\frac{mK^{1/2}}{N}\right)^{1/2} + \frac{mK^{1/2}}{N^{1/2}}\cdot \left(\frac{mK^{1/2}}{N}\right)^{-1/2}\right)
\ll mK^{1/2},
\end{split}
\end{equation}
where we have used \eqref{note} and \eqref{G''asymp}.

Next, we assume that $n\in J''_m\subseteq J_m$, where 
\begin{equation} \label{J''}
J''_m:=J_m\setminus J'(m)=\left[4c^{\ast}mr+\frac{c^{\ast}m\sqrt{t}}{\sqrt{2\pi hk\Delta_0N'}},4c^{\ast}mr+\frac{2c^{\ast}m\sqrt{t}N'}{\sqrt{2\pi hk\Delta_0}N^{3/2}}\right]
\end{equation}
and note that in this case 
\begin{equation} \label{xnote}
x\asymp N^{1/2} \quad \mbox{ if } A_m(n)\le x\le B_m(n).
\end{equation} 
If $n\in J''_m\subseteq J_m$, an estimate of the form $|G''_{m,n}(x)| \asymp \lambda$ doesn't hold on the entire $x$-interval. To make Lemma \ref{C} applicable,
we split this interval into subintervals such that an estimate of this form holds on each of them, except for one subinterval
which is so short that it can be treated trivially. To this end, 
we write
$$
I_m(n):= \left.\left(\frac{\left(4c^{\ast}\right)^{3/2}m t^{1/2}}{(2\pi |d|hk\Delta_0)^{1/2}\left(n-4mc^{\ast}r\right)}-1,\frac{\left(4c^{\ast}\right)^{3/2}m t^{1/2}}{(2\pi |d|hk\Delta_0)^{1/2}\left(n-4mc^{\ast}r\right)}\right.\right]
$$
and 
$$
I_m(n,\delta):= \left.\left(\frac{\left(4c^{\ast}\right)^{3/2}m t^{1/2}}{(2\pi |d|hk\Delta_0)^{1/2}\left(n-4mc^{\ast}r\right)}-2\delta,\frac{\left(4c^{\ast}\right)^{3/2}m t^{1/2}}{(2\pi |d|hk\Delta_0)^{1/2}\left(n-4mc^{\ast}r\right)}-\delta\right.\right]
$$
and set
$$
S_m(n):=I_m(n)\cap [A_m(n),B_m(n)] \quad \mbox{and} \quad S_m(n,\delta):=I_m(n,\delta)\cap [A_m(n),B_m(n)].
$$
We observe that the interval $[A_m(n),B_m(n)]$ can be split into $O\left(\log T\right)$ intervals of the form $S_m(n)$ or $S_m(n,\delta)$, where 
\begin{equation} \label{deltaest}
1\ll \delta\ll \sqrt{N}.
\end{equation}  

Estimating trivially gives 
\begin{equation} \label{triv}
\sum\limits_{n\in J_m''} \sum\limits_{x\in S_m(n)}  \frac{e\left(G_{m,n}(x)\right)}{\sqrt{2m|F_x''(y_{x,m,n})|}} \ll m^{1/2}K^{1/4},
\end{equation}
where we have used \eqref{Fprime} and \eqref{note}.
If $x\in S_m(n,\delta)$, then we compute using \eqref{J''}, \eqref{xnote} and the mean value theorem that 
\begin{equation*}
-(n-4c^{\ast}mr)^{2/3}|d|^{1/3}x^{2/3}+4c^{\ast}\cdot \frac{m^{2/3}t^{1/3}}{(2\pi hk\Delta_0)^{1/3}} \asymp \delta \cdot \left(\frac{m\sqrt{t}}{\sqrt{hk\Delta_0}N^{1/2}}\right)^{2/3} x^{-1/3} 
\asymp \frac{\delta m^{2/3}K^{1/3}}{N^{1/2}}.
\end{equation*}
Hence, using \eqref{preest}, it follows that
\begin{equation} \label{G''asymp2}
G_{m,n}''(x) \asymp \frac{mK^{1/2}}{\delta^{1/2}N^{3/4}}.
\end{equation}
Again using partial summation to remove the factor $1/\sqrt{F_x''(y_{x,m,n})}$ and Lemma \ref{C}, we deduce that
\begin{equation} \label{from}
\begin{split}
& \sum\limits_{n\in J_m''} \sum\limits_{x\in S_m(n,\delta)}  \frac{e\left(G_{m,n}(x)\right)}{\sqrt{2m|F_x''(y_{x,m,n})|}} \\ \ll & \frac{N^{1/2}}{m^{1/2}K^{1/4}}\cdot \left(\frac{\delta mK^{1/2}}{N^{1/2}} \cdot \left(\frac{mK^{1/2}}{\delta^{1/2}N^{3/4}}\right)^{1/2} + \frac{mK^{1/2}}{N^{1/2}} \cdot \left(\frac{mK^{1/2}}{\delta^{1/2}N^{3/4}}\right)^{-1/2}\right)
\ll mK^{1/2},
\end{split}
\end{equation}
where we have employed $|S_m(n,\delta)|\le \delta$, \eqref{note}, \eqref{deltaest} and \eqref{G''asymp2}. From \eqref{triv} and \eqref{from}, it follows that
\begin{equation} \label{follows}
\sum\limits_{n\in J_m''} \sum\limits_{A_m(n)\le x\le B_m(n)}  \frac{e\left(G_{m,n}(x)\right)}{\sqrt{2m|F_x''(y_{x,m,n})|}} \ll mK^{1/2}\log T.
\end{equation}

Now using $J_m=J_m'\cup J_m''$, \eqref{aftervandercorput} and \eqref{follows}, we get
\begin{equation} \label{get}
 \sum\limits_{n\in J_m} \sum\limits_{A_m(n)\le x\le B_m(n)}  \frac{e\left(G_{m,n}(x)\right)}{\sqrt{2m|F_x''(y_{x,m,n})|}}  \ll mK^{1/2}\log T.
\end{equation}
Combining \eqref{cs}, \eqref{rearrange} and \eqref{get}, we deduce that 
\begin{equation} \label{nontrivial}
\begin{split}
& \left|\sum\limits_{x\in J} \sum\limits_{y\in I(x)} e\left(Q^{\ast}(x,y)\cdot r+\frac{t}{\pi}\cdot \phi\left(\frac{\pi Q^{\ast}(x,y)}{2hk\Delta_0 t}\right)\right)\right|^2\\
= & 
O\left(MN^{3/2}+\frac{N^2}{M}+\frac{MN^{5/2}K^{3/2}}{T^2}+M^3N^{1/2}K^{1/2}+\frac{N^2}{M^{1/2}K^{1/4}}+MNK^{1/2}\log T\right).
\end{split}
\end{equation}
Choosing $M:=\left[N^{1/2}K^{-1/4}\right]$ to balance the second and last $O$-terms above, using $N\ll K$, and taking the square root, it follows that 
\begin{equation*}
\begin{split}
& K^{1/4}N^{-1/4}T^{-1/2}\left|\sum\limits_{x\in J} \sum\limits_{y\in I(x)} e\left(Q^{\ast}(x,y)\cdot r+\frac{t}{\pi}\cdot \phi\left(\frac{\pi Q^{\ast}(x,y)}{2hk\Delta_0 t}\right)\right)\right|\\
= & O\left(\frac{K^{7/8}}{T^{1/2}}\log T+\frac{K^{17/8}}{T^{3/2}}\right).
\end{split}
\end{equation*}
This is $O(T^{-\varepsilon})$, provided that $K:=T^{4/7-\varepsilon}$, which completes the proof. 

\medskip

\noindent
\textbf{Acknowledgement.} We express our gratitude to Professor Matti Jutila for pointing out some typos and suggesting some changes in an earlier version this paper. We thank Prof. A. Ivi\'c for encouragement   and for some useful suggestions. The third author would like to thank the Institute of Mathematical Sciences, Chennai for supporting her visits which enabled her to work on this project and would
like to mention that this paper constitutes part of her Ph D thesis work.

\end{document}